\documentclass[a4paper]{article}
\usepackage[russian]{babel}
\usepackage[T1]{fontenc}
\usepackage{pifont,setspace,hyperref,amsthm,amsmath,amssymb,color,multirow,graphicx,subfigure, color }
\usepackage{enumitem}
\usepackage{inputenc}
\usepackage{amsmath,latexsym,amsthm}
\usepackage{graphicx}
\usepackage{amssymb}
\begin{document}

\begin{center}
\textbf{\Large \textsc {On the positive fixed points of quartic operators}}\\[0.2 cm]
\textbf{Eshkabilov Yu.~Kh.~\footnote{Karshi state university,
Karshi, Uzbekistan. E-mail: yusup62@mail.ru}} \textbf{Nodirov
Sh.~D.~\footnote{Karshi state university, Karshi, Uzbekistan.
E-mail: shoh0809@mail.ru}}
\end{center}

%Below title, abstract and Keywords of the article will be given. First in Uzbek with Latine letters (in case, if you do not know Uzbek, editors will do it for you). Next, the same information should be given in Russian (in case, if you do not know Uzbek, editors will do it for you).
%
%
{\small
\begin{center}
\begin{tabular}{p{9cm}}
\textcolor{blue}{О положительных неподвижных точках квартических операторов}\\
В статье рассмотрен квартик оператор  с положительными коэффициентами
на $\mathbb{R}^2$. Исследованы положительные неподвижные точки
квартического оператора. Доказаны теоремы о количестве положительных
неподвижных точек квартического оператора.\\
\noindent \underline{Ключевые слова:} Неподвижная точка; квартический
оператор; метод Феррари; правила Декарта; трансляционно-инвариантные
меры Гиббса; дерево Кэли; модел; потенциал.
\end{tabular}\end{center} }
\vspace{0.5 cm}


\makeatletter
\renewcommand{\@oddhead}{\vbox{\hfill
{\it Eshkabilov Yu.~Kh.,~Nodirov Sh.~D.~ On the positive fixed
points of quartic operators}\hfill \thepage \hrule}}
\renewcommand{\@evenhead}{\vbox{\hfill
{ Bulletin of the Institute of Mathematics, 2019, \No 1,
ISSN-2181-9483}\hfill \thepage \hrule}}

\makeatother
% Label for the first page
%
\label{firstpage}

\section*{Introduction}

    It is well known that fixed point theory is a very important tool
for solving problems in Nonlinear Functional Analysis and as well as
to various theoretical and applied fields such as variational and
linear inequalities, the dynamic systems theory, nonlinear analysis
\cite{ap}. Nonlinear operators are connected with problems in
statistical physics, biology, thermodynamics, statistical mechanics
and so on \cite{ehr1}, \cite{to},\cite{cso}. One of the central
problem in statistical physics is the description of the set of
Gibbs measures for given models on Cayley trees.
Translation-invariant Gibbs measures for models with a continuum set
of spin values is reduced to find positive fixed points of nonlinear
integral operators \cite{ehr1}-\cite{ehr}. Shown that in
\cite{18},\cite{21} to find positive fixed point of nonlinear
integral operators require to find positive fixed points of
nonlinear operators on the finite dimensional space
$\mathbb{R}^{2}$.

In \cite{18} number of the positive fixed points of quadratic
operator on $\mathbb{R}^{2}$ was studied. A number of positive fixed
points of cubic operator on $\mathbb{R}^{2}$ were studied in
\cite{21}. In \cite{18},\cite{21} proven theorems, that
 quadratic and cubic operators
can have  one, two and  three  positive fixed points on
$\mathbb{R}^{2}$.

In this paper, we consider the   operator of degree fourth (quartic)
on $\mathbb{R}^{2}$.  We show that the existence of positive fixed
point of quartic operators and they can have up to five positive
fixed points. Besides, obtained the theorems for the number of
positive fixed points of quartic operators. Furthermore, gives
sufficient conditions for the number of positive fixed points of
quartic operators (table 2) and at the end gives examples to satisfy
conditions. Investigated quartic operators had arisen when examining
translation-invariant Gibbs measures for models in \cite{18} on the
Cayley tree of order $k=4$. In the paragraph 5 given model, which the
study translation-invariant Gibbs measure reduced the study the
fixed points of a quartic operator.

\section{Quartic operators on $\mathbb{R}_+^{2}$}

We introduce
$$\mathbb{R}_{+}^{2}=\{(x,y)\in \mathbb{R}^{2}: \,\ x\geq0, y\geq0\},\,\,\ \mathbb{R}_{>}^{2}=\{(x,y)\in \mathbb{R}^{2}: \,\ x>0, y>0\}.$$
We consider the following  operator  $\mathcal{Q}$ on the cone
$\mathbb{R}_+^2$:
\begin{equation}\label{i4.1}\mathcal{Q}(x,y) =
\left(\sum_{i=0}^4C_{4}^ia_{i}x^{4-i}y^i, \,\
\sum_{i=0}^4C_{4}^ib_{i}x^{4-i}y^i\right),
\end{equation}

where $C_{4}^i=\frac{4!}{(4-i)!i!}$ (binomial coefficient), $a_{i}>0$ and $b_{i}>0$ for all $i\in\{0,1,2,3,4\}$. This
operator is composed of two polynomials  of the fourth degree in two
variables, therefore we call the operator $\mathcal{Q}$ (\ref{i4.1})
a quartic. Clearly, an arbitrary non-trivial positive fixed point
$(x_0,y_0)\in\mathbb{R}_{+}^{2}$ of the quartic operator (QO)
$\mathcal{Q}$ is a strictly positive, i.e. $(x_0, y_0)\in
\mathbb{R}_{>}^{2} $. Denote by $N_{>}^{fix}(\mathcal{Q})$ the number of
 fixed points of the QO $\mathcal{Q}$ that belongs to
$\mathbb{R}_{>}^{2}$.

\textbf{Lemma 1.} \emph{If $\omega=(x_0,y_0)\in \mathbb{R}_{>}^{2}$
is a fixed point of the QO $\mathcal{Q}$, then $\xi_0=\frac{y_0}{x_0}$ is root of the
algebraic equation
\begin{equation}\label{e4.1}
a_{4}\xi^5+\left(4a_{3}-b_{4}\right)\xi^4+\left(6a_{2}-4b_{3}\right)\xi^3+(4a_{1}-6b_{2})\xi^2-(a_{0}-4b_{1})\xi-b_{0}=0.
\end{equation}
}

\medskip

{\sf Proof. }
 Let the point $\omega=(x_0,y_0)\in
\mathbb{R}_{>}^{2}$ be a fixed point of QO $\mathcal{Q}$. Then
$$\sum_{i=0}^4C_{4}^ia_{i}x_0^{4-i}y_0^i=x_0, \,\
\sum_{i=0}^4C_{4}^ib_{i}x_0^{4-i}y_0^i=y_0.$$

  From the
notation $\frac{y_{0}}{x_0} = \xi_0$ we obtain

$$x_0^{4}\left(\sum_{i=0}^4C_{4}^ia_{i}\xi_0^i\right)=x_0, \,\
x_0^{4}\left(\sum_{i=0}^4C_{4}^ib_{i}\xi_0^i\right)=\xi_{0}x_{0}.$$

Hence, we have
$$ \frac{1}{\xi_0}=\frac{\sum\limits_{i=0}^4C_{4}^ia_{i}\xi_0^i}{\sum\limits_{i=0}^4C_{4}^ib_{i}\xi_0^i}.$$

By the last equality we get
$$ a_{4}\xi_0^5+\left(4a_{3}-b_{4}\right)\xi_0^4+\left(6a_{2}-4b_{3}\right)\xi_0^3+(4a_{1}-6b_{2})\xi_0^2-(a_{0}-4b_{1})\xi_0-b_{0}=0. $$

Lemma 1 is proved. \hfill $\Box$\\

\textbf{Lemma 2.} \emph{ If $\xi_0$ is a positive root of the algebraic
equation (\ref{e4.1}), then the point $\omega_0=(x_0,\xi_{0}x_0)\in
\mathbb{R}_{>}^{2}$  and $\xi_0$ is a fixed point of the QO $\mathcal{Q}$, where
\begin{equation}\label{e4.2}
x_{0}=\frac{1}{\sqrt[3]{\sum\limits_{i=0}^4C_{4}^ia_{i}\xi_0^i}}.
\end{equation}
}

\medskip

{\sf Proof. }
 Let $\xi_0>0$ and $\xi_0$ is a root of the
equation (\ref{e4.1}). Put $y_0=\xi_0x_0$, where $x_0$ is a given by
the equality (\ref{e4.2}) and $\omega_0=(x_0,\xi_{0}x_0)$. From the
equality $y_0=\xi_0x_0$ we have

$$\sum_{i=0}^4C_{4}^ia_{i}x_0^{4-i}y_0^i=\sum_{i=0}^4C_{4}^ia_{i}x_0^{4-i}(\xi_0x_0)^i=x_0^4\left(\sum_{i=0}^4C_{4}^ia_{i}\xi_0^i\right)=
\frac{1}{\sqrt[3]{\sum\limits_{i=0}^4C_{4}^ia_{i}\xi_0^i}},$$ i.e.
$$\sum_{i=0}^4C_{4}^ia_{i}x_0^{4-i}y_0^i=x_0.$$
On the other hand
$$a_{4}\xi_0^5+\left(4a_{3}-b_{4}\right)\xi_0^4+\left(6a_{2}-4b_{3}\right)\xi_0^3+(4a_{1}-6b_{2})\xi_0^2-(a_{0}-4b_{1})\xi_0-b_{0}=0.$$
Then we get
$$
b_{0}+4b_{1}\xi_0+6b_{2}\xi_0^2+4b_{3}\xi_0^3+b_{4}\xi_0^4=a_{0}\xi_0+4a_{1}\xi_0^2+6a_{2}\xi_0^3+4a_{3}\xi_0^4+a_{4}\xi_0^5=$$
$$=\xi_0\left(a_{0}+4a_{1}\xi_0+6a_{2}\xi_0^2+4a_{3}\xi_0^3+a_{4}\xi_0^4\right)$$
From the last equality we have
$$\frac{\xi_0}{\sqrt[3]{\sum\limits_{i=0}^4C_{4}^ia_{i}\xi_0^i}}=
\frac{\sum\limits_{i=0}^4C_{4}^ib_{i}\xi_0^i}{\left(\sqrt[3]{\sum\limits_{i=0}^4C_{4}^ia_{i}\xi_0^i}\right)^{4}}=x_0^4\left(\sum\limits_{i=0}^4C_{4}^ib_{i}\xi_0^i\right)=\sum\limits_{i=0}^4C_{4}^ib_{i}x_0^{4-i}y_0^i=y_0.$$

 \hfill $\Box$\\

We put
$$\mu_0=a_{4},\,\ \mu_1=4a_{3}-b_{4},\,\,\ \mu_2=6a_{2}-4b_{3},\,\,\ \mu_3=4a_{1}-6b_{2},\,\,\,\,\ \mu_4=a_{0}-4b_{1},\,\,\ \mu_5=b_{0}$$
and define the polynomial $P_5(\xi)$ of degree 5:

\begin{equation}\label{add} P_5(\xi)=\mu_0\xi^5+\mu_1\xi^4+\mu_2\xi^3+\mu_3\xi^2+\mu_4\xi-\mu_5. \end{equation}

By the Lemmas 1, 2 the following corollary is correct.

\textbf{Corollary 1.}  \emph{The number of positive fixed points of the QO
$\mathcal{Q}$ equal to the number of positive roots of the polynomial
$P_5(\xi)$.}\\

\textbf{Lemma 3.}\emph{  The
 polynomial $P_5(\xi)$  has at least one positive root.
}

\medskip

{\sf Proof. }  Clearly, that $P_5(0)=-\mu_5<0$ and
$P_5(+\infty)=+\infty$. It means that there exists $c>0$ such that
$P_4(c)=0.$

 \hfill $\Box$

\textbf{Proposition 1.} \emph{ The polynomial $P_5(\xi)$ can have up
to five  positive roots. }

\medskip

{\sf Proof. } We have the following table for the number of
sign changes in the sequence of the coefficients of the polynomial $P_5(\xi)$ (Tab.
1). By using this table and the Descartes rule, we can conclude that
 the polynomial $P_5(\xi)$ can have up to five the positive roots  (see \cite{17}, pp. 27-29 ).

\bigskip
\begin{center}
\begin{tabular}{|c|c|c|c|c|c|c|c|}
  \hline
 $ P_5(\xi)$ & $\mu_0$& $\mu_1$ & $\mu_2$ & $\mu_3$ & $\mu_4$ &$\mu_5$& the number of sign changes \\
  \hline
  1. & + & - & + & - & + & - &5 \\

  \hline
\end{tabular}\end{center}
\begin{center}Table 1.\end{center}
\bigskip
\hfill $\Box$\\

\section{Sufficient conditions for number of positive fixed points of the QO $\mathcal{Q}$}

In this section we study a number of positive fixed points of the
quartic operator $\mathcal{Q}$.

We use the following designations
$$Q=\left(\frac{a}{3}\right)^3+\left(\frac{b}{2}\right)^2,$$
where
$$a=-\frac{p^2}{12}-r,\,\,\,\ b=-\frac{p^3}{108}+\frac{pr}{3}-\frac{q^2}{8},$$
and
$$p=\frac{15\mu_0\mu_2-6\mu_1^2}{25\mu_0^2},\,\,\,\ q=\frac{50\mu_0^2\mu_3+8\mu_1^3-30\mu_0\mu_1\mu_2}{125\mu_0^3
},$$
$$r=\frac{15\mu_0\mu_1^2\mu_2-50\mu_0^2\mu_1\mu_3-3\mu_1^4-125\mu_0^3\mu_4}{625\mu_0^4}.$$

When  $Q<0$, we define following number
$$z_0=2\sqrt{-\frac{a}{3}}\cos\left(\frac{2\pi}{3}+\frac{\alpha}{3}\right)-\frac{p}{3},$$
where
$$\cos\alpha=-\frac{b}{2}\left(-\frac{3}{a}\right)^{\frac{3}{2}},\,\,\,\
\alpha\in[0,\pi].$$

We also introduce
$$\xi^{ext}_{1,2}=\frac{1}{2}\left(\sqrt{2z_0}\pm\sqrt{2z_0-4\left(\frac{p}{2}+z_0+\frac{q}{2\sqrt{2z_0}}\right)}\right)-\frac{\mu_1}{5\mu_0},$$
$$\xi^{ext}_{3,4}=\frac{1}{2}\left(-\sqrt{2z_0}\pm\sqrt{2z_0-4\left(\frac{p}{2}+z_0-\frac{q}{2\sqrt{2z_0}}\right)}\right)-\frac{\mu_1}{5\mu_0}.$$
 Let us note that if $z_0>0$ then the numbers
 $\xi^{ext}_1,\xi^{ext}_2,\xi^{ext}_3,\xi^{ext}_4$ are real.

For $z_0>0$, denote $$\xi_{min}=\min\{\xi^{ext}_1, \xi^{ext}_2,
\xi^{ext}_3, \xi^{ext}_4\},\,\,\ \xi_{max}=\max\{\xi^{ext}_1,
\xi^{ext}_2, \xi^{ext}_3, \xi^{ext}_4\}.$$

\textbf{Theorem 1.} \emph{Let $Q<0$, $z_0>0$, $\xi_{max}<0$. Then
the QO $\mathcal{Q}$ has a unique fixed point in
$\mathbb{R}_{>}^{2}$, i.e., $N_{>}^{fix}(\mathcal{Q})=1$. }

\medskip

{\sf Proof. } By the Corollary 1 quite enough to find the positive
roots of the polynomial $P_5(\xi)$.\\
 Let $Q<0$. Then by the Cardano's formula
cubic equation $\eta^3+a\eta+b=0$ has three real roots. They have
the following form \cite{19}:
$$\eta_k=2\sqrt{-\frac{a}{3}}\cos\left(\frac{\alpha+2\pi(k-2)}{3}\right),\,\,\ k=\overline{1,2,3},$$
where
$$\cos\alpha=-\frac{b}{2}\left(-\frac{3}{a}\right)^{\frac{3}{2}},\,\,\,\
\alpha\in[0,\pi].$$ We have $\eta_3<\eta_1<\eta_2$. The number
$z_0=\eta_3-\frac{p}{3}(z_0>0)$ is the least root of the following
cubic equation:
\begin{equation}\label{cub}z^3+pz^2+\frac{p^2-4r}{4}z-\frac{q^2}{8}=0.\end{equation}
On the other hand if all roots are strictly positive of the equation
(\ref{cub}), then all roots $\omega_j,\,\ (j=\overline{1,..,4})$ of the
equation
$$\omega^4+p\omega^2+q\omega+r=0$$
 are real \cite{20}. By the Ferrari's method the numbers $\xi^{ext}_j=\omega_j-\frac{\mu_1}{5\mu_0}(j=\overline{1,..,4})$ are
roots of the equation $P_5'(\xi)=0$, i.e.
\begin{equation}\label{4.4}5\mu_0\xi^4+4\mu_1\xi^3+3\mu_2\xi^2+2\mu_3\xi+\mu_4=0.\end{equation}
Thus, the numbers $\xi^{ext}_j (j=\overline{1,..,4})$ are the local extreme points of the
function $P_5(\xi)$. Besides
$\xi_{max}=\max_{j\in\{1,...,4\}}\xi^{ext}_j<0$ and
$P_5(+\infty)=+\infty$. Therefore the function $P_5(\xi)$ is an
increasing on the set $(\xi_{max},+\infty)$. Consequently, by the
inequality $P_5(0)<0$ the polynomial $P_5(\xi)$ has a unique positive
root.

 \hfill $\Box$\\

\textbf{Theorem 2.} \emph{ Let $Q<0$, $z_0>0$, $\xi_{min}>0$. If the
polynomial $P_5(\xi)$ satisfies one of the following conditions}

$(a)\,\ P_5(\xi_{min})=0,$\\

$(b)\,\ P_5(\xi_{max})=0,$\\
 \emph{then the QO $\mathcal{Q}$ has at least two
fixed points in $\mathbb{R}_{>}^{2}$, i.e.,
$N_{>}^{fix}(\mathcal{Q})\geq2$. }

\medskip

{\sf Proof. } Let $Q<0$, $z_0>0$, $\xi_{min}>0$. Then all four roots
of the equation (\ref{4.4}) is real and they strictly positive. And so
the function $P_5(\xi)$ is an increasing on the set
$(-\infty,\xi_{min})\cup(\xi_{max},+\infty)$.

$(a)\,\ P_5(\xi_{min})=0.$ Then
$\max_{\xi\in(-\infty,\xi_{min}]}P_5(\xi)=P_5(\xi_{min})=0$ and the number
$\xi_1=\xi_{min}$ is the root of the polynomial $P_5(\xi)$. Besides
we know, that the function $P_5(\xi)$  is an increasing on the set
$(\xi_{max},+\infty)$. Therefore exist $\xi_2\in(\xi_{max},+\infty)$,
where $P_5(\xi_2)=0$. And so the polynomial $P_5(\xi)$ has at least
two strictly positive roots.

$(b)\,\ P_5(\xi_{max})=0.$ Then
$\min_{\xi\in[\xi_{max},+\infty)}P_5(\xi)=P_5(\xi_{max})=0$ and the number
$\xi_1=\xi_{max}$ is the root of the polynomial $P_5(\xi)$. The
function $P_5(\xi)$ is an increasing on the set
$[\xi_{max},+\infty)$, so exist $\xi_{k_0}\neq\xi_{max}$ which the
function $P_5(\xi)$ is decreasing on the set
$(\xi_{k_0},\xi_{max})$. Also we have $P_5(0)<0$
and $\max_{\xi\in[\xi_{k_0},\xi_{max}]}P_5(\xi)=P_5(\xi_{k_0})>0$. Therefore the polynomial
$P_5(\xi)$ has another unique positive root $\xi_2 $ on the set $(0,\xi_{k_0}).$

 \hfill $\Box$\\

\textbf{Theorem 3.} \emph{ Let $Q<0$, $z_0>0$, $\xi_{min}>0$. If for
the polynomial $P_5(\xi)$ following condition}

$(c)\,\ P_5(\xi_{min})>0,\,\,\ P_5(\xi_{max})<0$\\
\emph{is satisfies, then the QO $\mathcal{Q}$ has at least three
fixed points in $\mathbb{R}_{>}^{2}$, i.e.,
$N_{>}^{fix}(\mathcal{Q})\geq3$. }

\medskip

{\sf Proof. } Let $Q<0$, $z_0>0$, $\xi_{min}>0$. Then all four roots of the
equation (\ref{4.4}) is real and they strictly positive. And so the
function $P_5(\xi)$ is an increasing on the set
$(-\infty,\xi_{min})\cup(\xi_{max},+\infty)$.

$(c)$  Let $ \,\ P_5(\xi_{min})>0,\,\,\ P_5(\xi_{max})<0.$ We have known that
$P_5(0)<0$ and $P_5(+\infty)=+\infty.$ It indicates that  the
polynomial $P_5(\xi)$ has at least three positive roots, such that
$\xi_1\in(0,\xi_{min}), \,\ \xi_2\in(\xi_{min}, \xi_{max})$ and
$\xi_3\in(\xi_{max},+\infty).$

 \hfill $\Box$ \\

In the next, we assume that  $Q<0$, $z_0>0$. Then we can define
following set of real
 numbers (set of extremal points of the function $P_5(\xi)$ ):
$$\mathfrak{A}=\{\xi^{ext}_1,\xi^{ext}_2,\xi^{ext}_3,\xi^{ext}_4\}.$$

 Points of the set of
$\mathfrak{A}$ are real and they are different. Therefore  we  can have denote following
$$\lambda_1=\xi_{min},\,\,\ \lambda_2=\min \left(\mathfrak{A}\backslash\{\xi_{min}\}\right),\,\,\,\ \lambda_3=\max\left(\mathfrak{A}\backslash\{\xi_{max}\}\right),\,\,\ \lambda_4=\xi_{max}.$$

Consequently, we have the following relation
$$\lambda_1<\lambda_2<\lambda_3<\lambda_4.$$

It follows the function $P_5(\xi)$ is an increasing (decreasing) on
the set $(-\infty,\lambda_1)\cup
(\lambda_2,\lambda_3)\cup(\lambda_4, +\infty)$ $\bigl((\lambda_1,
\lambda_2)\cup(\lambda_3, \lambda_4)\bigr)$. The function
$P_5(\xi)$ has a local maximum value
 at the points $\lambda_1,\lambda_3$ and local minimum values at the points $\lambda_2$ and
 $\lambda_4$. Using these property and by the corollary 1, we can give the following
 table:\\

\bigskip
\begin{center}
\begin{tabular}{|c|c|c|}
  \hline
\multicolumn{3}{|c|}{ On the number of the positive fixed points of
the quartic
operators, when $Q<0$, $z_0>0$, $\lambda_1>0.$ }\\
\hline
   & Sufficient condition:  & $N_{>}^{fix}(\mathcal{Q})$ \\
  \hline
1.&$\,\ P_5(\lambda_1)<0,\,\,\ P_5(\lambda_3)<0$&1 \\
\cline{1-2} \cline{1-2} \hline
 2.&$\,\ P_5(\lambda_1)=0,\,\,\ P_5(\lambda_3)<0$& \\
\cline{1-2}

3.&$\,\ P_5(\lambda_1)<0,\,\,\ P_5(\lambda_3)=0$&\\ \cline{1-2}

4. &$\,\ P_5(\lambda_2)>0,\,\,\ P_5(\lambda_4)=0$&2\\
\cline{1-2}

5.&$\,\ P_5(\lambda_2)=0,\,\,\ P_5(\lambda_4 )>0$&\\ \cline{1-2}

 6.&$\,\ P_5(\lambda_2)>0,\,\,\ P_5(\lambda_4)=0$ &  \\ \cline{1-2}
  \hline
 7.&$\,\ P_5(\lambda_1)=0,\,\,\ P_5(\lambda_4)=0$ & \\ \cline{1-2}

8.&$\,\ P_5(\lambda_1)=0,\,\,\ P_5(\lambda_3)=0$ & \\
\cline{1-2}

9.&$ \,\ P_5(\lambda_2)=0,\,\,\ P_5(\lambda_4)=0$ & \\
\cline{1-2}

10.&$\,\ P_5(\lambda_2)>0,\,\,\ P_5(\lambda_4)<0$ & 3\\
\cline{1-2}

11&$\,\ P_5(\lambda_1 )<0,\,\,\ P_5(\lambda_3)>0, \,\,\
P_5(\lambda_4)<0$ &  \\ \cline{1-2}

12.&$\,\ P_5(\lambda_1)>0,\,\,\ P_5(\lambda_2)<0, \,\,\
P_5(\lambda_4)>0$ &   \\ \cline{1-2}
   \hline
13.&$\,\ P_5(\lambda_1)=0,\,\,\ P_5(\lambda_3)>0,\,\,\
P_5(\lambda_4)<0$ &\\ \cline{1-2}

14.&$\,\ P_5(\lambda_1)>0,\,\,\ P_5(\lambda_3)=0$ & \\
\cline{1-2}

15.&$\,\ P_5(\lambda_2)=0,\,\,\ P_5(\lambda_4)<0$ & 4\\
\cline{1-2}

16.&$\,\ P_5(\lambda_1)>0,\,\,\ P_5(\lambda_2)<0,\,\,\
P_5(\lambda_4)=0$ &  \\ \cline{1-2}
   \hline
17.&$\,\ P_5(\lambda_1)>0,\,\,\ P_5(\lambda_2)<0,\,\,\
P_5(\lambda_3)>0,\,\,\ P_5(\lambda_4)<0$&  5 \\ \cline{1-2}
  \hline
\end{tabular}\\
\end{center}
\begin{center}Table 2.\end{center}
\bigskip

\section{ Examples}
In this section we give the concrete examples, which the QO
$\mathcal{Q}$ has the three and the five positive fixed points.\\

\textsc{Example 1.}  We consider an example for QO $\mathcal{Q}$,
which satisfies the 7th condition  in the table 2. We define QO
$\mathcal{Q}$ by the equality:
\begin{equation}\label{e1}
\mathcal{Q}(x,y)=\left(52x^4+8x^3y+42x^2y^2+4xy^3+y^4,\,\
16x^4+4x^3y+72x^2y^2+5xy^3+14y^4\right).\end{equation}
 Then for the QO $\mathcal{Q}$ we have
$$a_0=52,\,\ a_1=2,\,\ a_2=7,\,\ a_3=1,\,\ a_4=1,\,\ b_0=16,\,\ b_1=1,\,\ b_2=12,\,\ b_3=\frac{5}{4},\,\ b_4=1.$$
Therefore we get
$$\mu_0=a_{4}=1,\,\ \mu_1=4a_{3}-b_{4}=-10,\,\,\ \mu_2=6a_{2}-4b_{3}=37,\,\ \mu_3=4a_{1}-6b_{2}=-64,$$
$$ \mu_4=a_{0}-4b_{1}=52,\,\
\mu_5=b_{0}=16.$$
 For the QO $\mathcal{Q}$ (\ref{e1}) the polynomial
${{P}_{5}}(\xi)$ has the following form
$$P_5(\xi)=\xi^5-10\xi^4+37\xi^3-64\xi^2+52\xi-16.$$

Then we have
$$ Q=-\frac{7}{12500},\,\,\,\ z_0=\frac{1}{20}\left(13-\sqrt{105}\right ),$$

$$\xi_1^{ext}=1,\,\,\ \xi_2^{ext}=2,\,\,\ \xi_3^{ext}=\frac{1}{10}\left(25 - \sqrt{105}\right),\,\,\  \xi_4^{ext}=\frac{1}{10}\left(25 + \sqrt{105}\right).$$

By the denotes we have
$$\lambda_1=1,\,\,\ \lambda_2=\frac{1}{10}\left(25 - \sqrt{105}\right),\,\,\ \lambda_3=2,\,\,\  \lambda_4=\frac{1}{10}\left(25 + \sqrt{105}\right).$$
So
 $Q<0,\,\ z_0>0,\,\ \lambda_1>0,\,\,\ P_5(\lambda_1)=P_5(1)=0,\,\,\ P_5(\lambda_3)=P_5(2)=0.$ By the 7th condition  in the table 2, the  QO
 $\mathcal{Q}$ has three positive fixed points.

  It is easy to verify that the polynomial $P_5(\xi)$ has three positive roots and they have the following form:
$$\xi_{1,2}=1,\,\,\ \xi_{3,4}=2\,\,\ \xi_5=4.$$

 By the Corollary 1 the QO $\mathcal{Q}$ (\ref{e1})  has three
positive fixed points.\\

\textsc{Example 2.}  Now, we consider  the following QO $\mathcal{Q}$:

\begin{equation}\label{e2}\mathcal{Q}(x,y)=\left(31x^4+\frac{1}{2}x^3y+\frac{43}{3}x^2y^2+\frac{1}{4}xy^3+\frac{1}{5}y^4,
10x^4+x^3y+31x^2y^2+\frac{2}{3}xy^3+3y^4\right).\end{equation}

Then

$$P_5(\xi)=\frac{1}{5}\xi^5-\frac{11}{40}\xi^4+\frac{41}{3}\xi^3-\frac{61}{2}\xi^2+30\xi-106.$$

Also $$a=-\frac{7}{3},\,\ b=-\frac{20}{27},\,\ Q=-\frac{1}{3}, $$
and $$z_0=\frac{1}{8}.$$

Consequently,
$$\lambda_1=1,\,\,\ \lambda_2=2,\,\,\ \lambda_3=3,\,\,\  \lambda_4=5.$$
The following relations are true :
$$Q=-\frac{1}{3}<0,\,\,\ z_0=\frac{1}{8}>0,\,\,\ \lambda_1=1>0, $$
$$P_5(\lambda_1)=P_5(1)=\frac{37}{60}>0,\,\,\ P_5(\lambda_2)=P_5(2)=-\frac{4}{15}<0,$$
$$ P_5(\lambda_3)=P_5(3)=\frac{7}{20}>0,\,\,\ P_5(\lambda_4)=P_5(5)=-\frac{95}{12}<0,$$

By the 17th condition  in the table 2, the  QO
 $\mathcal{Q}$ (\ref{e2}) has five positive fixed points.

\section{Application}

In this section, we show that the number of fixed points of the operator $ \mathcal{Q} $ (\ref{i4.1}) corresponds to a non-trivial positive fixed point of the Hammerstein integral operator, which plays an important role in the theory of Gibbs measures.

A Cayley tree $\Gamma^k=(V,L)$ of order $k\geq 1$ is an infinite
homogeneous tree, i.e. a graph without cycles, with exactly $k+1$
edges incident to each vertices. Here $V$ is the set of vertices and
$L$ that of edges.

Consider models where the spin takes values in the set $[0,1]$, and
is assigned to the vertices of the tree. For $A\subset V$ a
configuration $\sigma_A$ on $A$ is an arbitrary function
$\sigma_A:A\to [0,1]$. Denote $\Omega_A=[0,1]^A$ is the set of all
configurations on $A$.

We consider the Hamiltonian of the model in \cite{21}:
\begin{equation}\label{m}
 H(\sigma)=-J\sum_{\langle x,y\rangle\in L}
\xi_{\sigma(x), \sigma(y)}, \,\ \sigma\in\Omega_{V}
\end{equation}
where $J \in R\setminus \{0\}$ and $\beta=\frac{1}{T},\,\ T>0$ is
temperature, $\xi: (u,v)\in [0,1]^2\to \xi_{uv}\in \mathbb{R}$ is a
given bounded, measurable function. As usually, $\langle x,y\rangle$
stands for the nearest neighbor vertices on the Cayley three of
order four.

We introduce
$$C_+[0,1]=\{f\in C[0,1]: f(x)\geq 0\}, \,\ C_+^0[0,1]=C_+[0,1]\setminus \{\theta\equiv0\},\,\ C_>[0,1]=\{f\in C[0,1]: f(x)>0\}.$$
Consider the following  operator \cite{ehr1} $R_{k}$ on $C_+[0,1]$ defined by:
$$(R_{k}f)(t)=\left({\int\limits_0^1K(t,u)f(u)du\over \int\limits_0^1
K(0,u)f(u)du}\right)^k, \,\
k\in\mathbb{N}$$\\
where $K(t,u)=\exp(J\beta \xi_{tu}), f(t)>0, t,u\in [0,1].$\\

In \cite{re}-\cite{ehr} was shown that a translation-invariant Gibbs measure of the model $H$ (\ref{m})  corresponds to a solution $f(t)\in
C_+^0[0,1]$  of the following equation:
\begin{equation}\label{e1.2}
(R_{k}f)(t)=f(t).
\end{equation}

For every $k\geq2$ we define the Hammerstein integral operator $H_k$ acting in $C_+[0, 1]$ as
following:
$$(H_k f)(t) = \int\limits_0^1 K(t, u) f^k(u)du.$$

\textbf{Lemma 4.} \emph{Let $k \geq 2$.
The equation (\ref{e1.2}) has a non-trivial positive solution iff the Hammerstein integral operator has a non-trivial positive fixed point and $$N_+^{fix}(R_k)=N_+^{fix}(H_k),$$
where $N_+^{fix}(T)$ is a number of non-trivial positive fixed points of the operator $T$.
}

\medskip

{\sf Proof. } At first, we notice that the equation (\ref{e1.2}) has at least one solution in $C_>[0,1]$ (see Theorem 3.5 in \cite{ehr}).\\
{\it Necessity.} Let $k\geq2$ and $f(t)\in C_+^0[0,1]$ be a solution of the  (\ref{e1.2}), then $f(0)=1$. Define the following function \begin{equation}\label{g1}g(t)=\frac{\sqrt[k]{f(t)}}{\sqrt[k-1]{\int\limits_0^1K(0,u)f(u)du}},\end{equation}
where $\int\limits_0^1K(0,u)f(u)du>0$. From $f(0)=1$, $g$ corresponds to exactly one function $f$.
Then we have  $$(H_kg)(t)=g(t).$$
{\it Sufficiency.}  Let $k\geq2$ and $g=g(t)\in C_+^0[0,1]$ be a fixed point of the operator $H_k$.  It is clear that  $\int\limits_0^1K(0,u)g^k(u)du=g(0)>0$.
Define non-trivial positive continuous function    \begin{equation}\label{f1}f(t)=\biggl(\frac{g(t)}{g(0)}\biggl)^k.\end{equation}
Because of $g$ is a fixed point of the operator $H_k$, the last equality is well-defined, i.e., each function $f$
corresponds to exactly one function $g$.

Then
$$\left(R_kf\right)(t)=\left(\frac{\int\limits_0^1K(t,u)f(u)du }{\int\limits_0^1
K(0,u)f(u)du}\right)^k=\left(\frac{\int\limits_0^1K(t,u)g^k(u)du}{\int\limits_0^1K(0,u)g^k(u)du}\right)^k=\left(\frac{g(t)}{g(0)}\right)^k=f(t).$$

It is easy to check that the relations (\ref {g1}) and (\ref{f1}) establish a one-to-one correspondence between
 $$F_\alpha=\{f(t)\in C_+^0[0,1]: R_kf=f\}\,\,\,\ \mbox{and}\,\,\,\ G_\alpha=\{g(t)\in  C_+^0[0,1]: H_kg=g\}.$$
 Consequently $N_+^{fix}(R_k)=N_+^{fix}(H_k).$

 This completes
the proof.
 \hfill $\Box$

 \vspace{0.2cm}

From Lemma 4 we can conclude that the number of the non-trivial positive fixed points of the operator $H_k$ is equal to the number of  translation-invariant Gibbs measures of the model $H$ (\ref{m}).

Let functions $\varphi_1(t), \,\ \varphi_2(t)$ and $\phi_1(t), \,\ \phi_2(t)$
belong to $C_>[0,1]$. We consider
the Hamiltonian $H$ (\ref{m}) on the Cayley tree of order $k=4$ with the
function of potential
\begin{equation}\label{p1}\xi_{t,u}=\frac{1}{J\beta}\ln\biggl(\phi_1(t)\varphi_1(u)+\phi_2(t)\varphi_2(u)\biggl).\end{equation}
 In a similar manner, we can prove that for the number non-trivial positive fixed points  of the integral operator $H_4$  the following equality holds, as in the \cite{21} (see Lemma 7):

$$N_{+}^{fix}(H_4)=N_{>}^{fix}(\mathcal{Q}).$$

By the last equality we can conclude that number of nontrival positive fixed points
of the quartic operator $\mathcal{Q}$ (\ref{i4.1})  on the $\mathbb{R}_>^2$ is equal
to number of translation-invariants Gibbs measures for the model
$H$ (\ref{m}) on the Cayley tree of order four.\\

{\it Open problem.} It has not  been constructed yet the exact function of potential $\xi_{t,u}$ (\ref{p1}) for satisfy the conditions of the Theorems 2-3. This requires constructing non-trivial positive continuous functions  $\varphi_1(t), \,\ \varphi_2(t)$ and $\phi_1(t), \,\ \phi_2(t)$ in such a way that the corresponding coefficients
$$a_i=\int\limits_{0}^{1}\varphi_{1}\left( u
\right)\phi _{1}^{k-i}\left( u \right)\phi _{2}^{i}\left( u
\right)du, \,\,\,\
b_{i}=\int\limits_{0}^{1}\varphi_{2}\left( u \right)\phi _{1}^{k-i}\left( u \right)\phi _{2}^{i}\left( u \right)du,\,\,\ i=\overline{0,..,4}$$
of the nonlinear operator $\mathcal{Q}$ need satisfy the conditions of the Theorems 2-3 (as in the above examples 1-2). \\

\vspace{0.1cm}
{\bf Acknowledgements.} We would like to thank reviewers for their insightful comments on the paper, as these comments led us to an improvement of the work.

%   label for the last page of the article
%
\label{lastpage}

\end{document}